\documentclass[opre,nonblindrev]{informs3wp}

\OneAndAHalfSpacedXI % current default line spacing
\usepackage{endnotes}
\let\footnote=\endnote

\usepackage{natbib}
 \bibpunct[, ]{(}{)}{,}{a}{}{,}%
 \def\bibfont{\small}%
\TheoremsNumberedThrough  % Preferred (Theorem 1, Lemma 1, Theorem 2)
\ECRepeatTheorems

\EquationsNumberedThrough % Default: (1), (2), ...
\usepackage[pdfpagemode={UseOutlines},bookmarks=true,bookmarksopen=true,
bookmarksopenlevel=0,bookmarksnumbered=true,hypertexnames=false,
colorlinks,linkcolor={blue},citecolor={black},urlcolor={red},
pdfstartview={FitV},unicode,breaklinks=true]{hyperref}
\hypersetup{urlcolor=blue, colorlinks=true}

\usepackage{multirow}% http://ctan.org/pkg/multirow
\usepackage{hhline}% http://ctan.org/pkg/hhline
\usepackage{footnote}
\usepackage[ruled,vlined]{algorithm2e}
\usepackage{algorithmic}
\newcommand{\comments}[1]{{\color{blue}\textit{$\#$ #1}}}

\newcommand{\transpose}{{\mbox{\tiny T}}}

\newcommand{\cV}{{\mathcal{V}}}

\newcommand{\cS}{{\mathcal{S}}}

\newcommand{\cO}{{\mathcal{O}}}

\newcommand{\bbe}{{\textbf{e}}}

\newcommand{\bB}{\textbf{B}}

\newcommand{\bY}{\textbf{Y}}
\newcommand{\bx}{\textbf{x}}
\newcommand{\by}{\textbf{y}}

\newcommand{\bU}{\textbf{U}}

\newcommand{\bH}{\textbf{H}}

\newcommand{\ba}{\textbf{a}}

\newcommand{\bbR}{\mathbb{R}}

\newtheorem{note}{Note}

\newcommand{\MMNL}{\textsc{\tiny MMNL}}
\newcommand{\GEV}{\textsc{\tiny GEV}}

\usepackage{mathtools}

\newif\ifnotes\notestrue
\def\ctien#1{{\color{red}{\textsc{[#1]}}}}
\def\htien#1{}

\begin{document}
%%%%%%%%%%%%%%%%

% Outcomment only when entries are known. Otherwise leave as is and
%default values will be used.
%\setcounter{page}{1}
%\VOLUME{00}%
%\NO{0}%
%\MONTH{Xxxxx}% (month or a similar seasonal id)
%\YEAR{0000}% e.g., 2005
%\FIRSTPAGE{000}%
%\LASTPAGE{000}%
%\SHORTYEAR{00}% shortened year (two-digit)
%\ISSUE{0000} %
%\LONGFIRSTPAGE{0001} %
%\DOI{10.1287/xxxx.0000.0000}%

% Author's names for the running heads
% Sample depending on the number of authors;
% \RUNAUTHOR{Jones}
% \RUNAUTHOR{Jones and Wilson}
% \RUNAUTHOR{Jones, Miller, and Wilson}
% \RUNAUTHOR{Jones et al.} % for four or more authors
% Enter authors following the given pattern:
\RUNAUTHOR{Dam T.T., Ta T.A and Mai T.}

% Title or shortened title suitable for running heads. Sample:
% \RUNTITLE{Bundling Information Goods of Decreasing Value}
% Enter the (shortened) title:
    \RUNTITLE{Maximum Capture Facility Location under Generalized Extreme Value Models}

% Full title. Sample:
% \TITLE{Bundling Information Goods of Decreasing Value}
% Enter the full title:
\TITLE{Submodularity and Local Search Approaches for Maximum Capture Problems under Generalized Extreme Value Models}

% Block of authors and their affiliations starts here:
% NOTE: Authors with same affiliation, if the order of authors allows,
%should be entered in ONE field, separated by a comma.
%\EMAIL field can be repeated if more than one author
\ARTICLEAUTHORS{%
\AUTHOR{Tien Thanh Dam}
\AFF{ORLab, Faculty of Computer Science, Phenikaa University, Yen Nghia, Ha Dong, Hanoi, VietNam \EMAIL{thanh.damtien@phenikaa-uni.edu.vn}} %, \URL{}}
\AUTHOR{Thuy Anh Ta}
\AFF{ORLab, Faculty of Computer Science, Phenikaa University, Yen Nghia, Ha Dong, Hanoi, VietNam \EMAIL{
anh.tathuy@phenikaa-uni.edu.vn}}
\AUTHOR{Tien Mai}
\AFF{School of Computing and Information Systems, Singapore Management University, \EMAIL{atmai@smu.edu.sg}}
% Enter all authors
} % end of the block

\ABSTRACT{%
We study  the maximum capture problem in facility location under random utility models, i.e., the problem of seeking to locate new facilities in a competitive market such that the captured user demand is maximized, assuming that each customer chooses among all available facilities according to a random
utility maximization model. We employ the generalized extreme value (GEV) family of discrete choice models and show that the objective function in this context is monotonic and submodular. This finding implies that a simple greed heuristic can  always guarantee an $(1-1/e)$ approximation solution.
We further develop a new algorithm combining a greedy heuristic, a gradient-based local search  and an exchanging procedure to efficiently solve the  problem.   
 We conduct experiments using instances of difference sizes and under different discrete choice models, and we show that our approach significantly outperforms prior approaches in terms of both returned objective value and CPU time.  
Our algorithm and theoretical findings can be applied to the maximum capture problems under various random utility models in the literature, including the popular multinomial logit, nested logit, cross nested logit, and the mixed logit models.
}%

% Sample
%\KEYWORDS{deterministic inventory theory; infinite linear programming duality;
%  existence of optimal policies; semi-Markov decision process; cyclic schedule}

% Fill in data. If unknown, outcomment the field
%\HISTORY{This paper wasfirst submitted on April 12, 1922 and has been with the authors for 83 years for 65 revisions.}

\maketitle

%%%%%%%%%%%%%%%%%%%%%%%%%%%%%%%%%%%%%%%%%%%%%%%%%%%%%%%%%%%%%%%%%%%%%%

% Samples of sectioning (and labeling) in OPRE
% NOTE: (1) \section and \subsection do NOT end with a period
% (2) \subsubsection and lower need end punctuation
% (3) capitalization is as shown (title style).
%
%\section{Introduction.}\label{intro} %%1.
%\subsection{Duality and the Classical EOQ Problem.}\label{class-EOQ} %% 1.1.
%\subsection{Outline.}\label{outline1} %% 1.2.
%\subsubsection{Cyclic Schedules for the General Deterministic SMDP.}
%  \label{cyclic-schedules} %% 1.2.1
%\section{Problem Description.}\label{problemdescription} %% 2.

% Text of your paper here

\section{Introduction}
In the last decade, the facility location problem in a competitive market has received a growing attention. In practice, modelling critical managerial decisions related to infrastructure planning, such as finding locations to locate new retail, service or product facilities in a market, often lead to facility location problems. The competitive facility location problem deals with a decision of  selecting locations to open new facilities in a market to maximize the captured demand of users, where a set of incumbent competitors are already operating in order.  There are two aspects that need to be considered in this problem, namely, the demand of customers and the competitors in the market. Customers are independent decision makers and their choices among different facilities might be based on a given utility that they assign to each location. Such utilities might be a function of facility attributes/features, e.g., distances, prices and transportation costs.

There are several ways to define and estimate customer demand \citep{Berman2009}. In this work, we focus on a probabilistic approach, i.e., customer demand is 
captured by a probability model that assigns choice probabilities to the facilities. The random utility maximization (RUM) framework \citep{McFa1973,Ben-Akiva1986} is convenient and popular  in the context. This framework is based on the assumption that each facility is associated with a random utility, which can be determined by the features/attributes of the facility. The RUM principle assumes that each customer selects a facility by maximizing his/her utilities. This way of modeling allows for predicting the probability that a customer selects a  facility. The facility location problem then  becomes the problem of locating new facilities in a competitive market to maximize an expected captured demand function, where customers selects a facility (a new facility or one from the competitors)  according to a RUM model. Thus,  the problem is also called as the \textit{maximum capture
problem} (MCP). 

To the best of our knowledge, existing related studies in the literature only employ the multinomial logit (MNL) or its mixed version (mixed logit model - MMNL) \citep{Benati2002,Haase2009,Haase2013}. However, it is well-known that the MNL retains the independence from irrelevant alternatives (IIA) property, 
which does not hold in many contexts \citep{McFa81,McFaTrai00}. On the other hand, the generalized  extreme value (GEV) family provides  flexible ways to relax the IIA property and capture the correlation between alternative utilities \citep{McFa81}. However, under the GEV family, most of the important properties that have been used to develop solution methods for the MCP under the MNL and MMNL models do not hold, or have not been proved to be true. More precisely, the objective function under the GEV family does not have a linear fractional structure, thus it is difficult to formulate the MCP into a mixed-integer linear program (MILP) as in prior work \citep{Benati2002,Zhang2012}. Moreover, under the GEV family, the objective function of the continuous relaxation is not either concave or convex, making the outer-approximation methods \citep{Ljubic2018outer,Mai2019Assortment} not applicable. Furthermore, since the structure of the objective function is driven by a GEV choice probability generating function, which may not have a closed form and could be complicated, it is not clear whether the objective function is submodular or not.  All the above remarks make the MCP under the GEV family challenging.
We tackle this challenge this in paper. 

Before presenting our contributions in detail, we note that, from now on, when saying a ``GEV model'', we refer to any choice model in the GEV family.
Each GEV model can be determined by a choice probability generating function (CPGF) $G(\cdot)$ \citep{Fosgerau2013}
(see our detailed definition in the next section).

\noindent\textbf{Our contributions:} 
In this paper, we formulate and solve the MCP under any GEV models.
 We leverage the properties of the CPGFs of GEV models \citep{McFa81,DalyBier06} and show that  the objective function in the context is  monotonic increasing  and submodular. These properties are already known for the MCP under MNL \citep{Benati2002} and now we show that they also hold for any GEV models. 
 The monotonicity  and submodularity also imply that the MCP subjecting to a cardinality constraint, even though being \textit{NP-hard}, always admits an $(1 - 1/e)$ approximation algorithm. In other words, a simple greedy heuristic always returns a solution whose value is  at least  $(1 - 1/e)$ ($\approx 0.632$) times the optimal values \citep{Nemhauser1978analysis}. 
 
 To further enhance the greedy heuristic (GH), we develop a new algorithm that adds a gradient-based local search and exchanging procedures  to the GH. While the latter is simply based on steps of  exchanging a location in a set of chosen locations with one outside of the set to get a better objective value, the former is motivated by the fact that if we formulate the MCP as a binary program, then the objective function is differentiable and we can make use of gradient information to direct the search. The gradient-base location search is an iterative procedure in which at each iteration  we solve a subproblem to (hopefully) find a better candidate solution,  and we show that such subproblems are solvable  in polynomial-time. Our algorithm can be used to solve problems under any GEV models and  under the MMNL model.

We conduct experiments using some datasets from the recent literature, including  real-life large-scale instances from a \textit{park-and-ride} location problem in New York City \citep{Holguin2012new}. We compare our algorithm, named as GGX (stands for \textbf{G}reedy Heuristic, \textbf{G}radient-based Local Search, and \textbf{Ex}changing) with some state-of-the art approaches from recent literature, i.e., the Branch \& Cut method proposed by \cite{Ljubic2018outer} and  outer-approximation algorithms \citep{Bonami2011_BB_MIP,Mai2019Assortment}. 
Experiments based on MNL, MMNL and nested logit instances  show that our algorithm remarkably outperforms the other approaches, in terms of both returned objective value and CPU time.

\noindent\textbf{Literature review:}  
The GEV family \citep{McFa81} covers  most of the discrete choice models in the demand modeling and operations research literature. Among existing GEV models, the MNL is the simplest and most popular one. It is also well-known that the MNL retains the  IIA property, which implies that the ratio between the choice probabilities of two facilities will not change no matter what other facilities are available or what attributes that other facilities have. This property has been considered as a limitation of the MNL model and should be relaxed in many applications \citep{McFaTrai00}. There are several GEV models that relax this property and provide flexible ways to model the correlation between  alternatives. For example, the nested logit   \citep{BenA73}, the cross-nested logit  \cite{VovsBekh98}, the generalized nested logit \citep{WenKopp01}, the paired combinatorial logit \citep{KoppWen00}, the ordered generalized extreme value \cite{Smal87}, the specialized compound generalized extreme value models \citep{WheBatFowDal02} and network-based GEV \citep{DalyBier06,MaiFreFosBas15_DynMEV} models. GEV models, in particular
the cross-nested and network GEV models, are fully flexible, in the sense that these models can approximate any random utility maximization models \citep{Fosgerau2013}. 
Beside the GEV family,the MMNL is also an alternative to relax the IIA property. This model  extends the MNL by assuming that choice parameters are random.
Similar to  GEV models, the MMNL is also able to approximate any
random utilities choice model \citep{McFaTrai00}. However, the choice probabilities given by the MMNL model have no closed form and often require simulation to approximate. Thus,  the estimation and the application of this model is expensive
in many contexts. 
 
In the context of the MCP, most  existing studies  focus on the MNL model due to its simplicity. \cite{Benati2002} seem the first to introduce the MCP under the MNL model. They propose three methods to compute upper bounds along with a branch-and-bound method to solve small instances. The first method is based on the concavity of the continuous relaxation of the objective function. They show the  submodularity of the objective function and use this property to develop the second method. The third method is an equivalent mixed-integer linear program (MILP),  which is
based on the fact that the objective function has a linear fractional structure  and can be linearized using additional additional variables. 
 \cite{Benati2002} also introduced a simple variable neighborhood search (VNS) method to solve instances with more than 50 potential locations. Some alternative MILP models, afterwards, have been proposed by \cite{Haase2009} and \cite{Zhang2012}. \cite{Haase2013} give an  evaluation and comparison between the proposed MILP models and conclude that the MILP model from \cite{Haase2009} is the most efficient one. \cite{Freire2015} strengthen the MILP reformulation of \cite{Haase2009} by using some tighter coefficients in some inequalities and also propose a new branch-and-bound algorithm to deal with the problem. 
More recently, \cite{Ljubic2017} propose a branch-and-cut method that combines two types of cutting planes, namely, outer-approximation (OA) cuts and submodular cuts. The first type of cuts is relied on the fact that the objective function of the continuous relaxation of the problem is concave and differentiable and the second type is based on the  submodularity and separability properties of objective function. Their branch-and-cut method is an iterative procedure where cuts are generated for every demand points and a linear programming (LP) relaxation is solved at each iteration. \cite{Mai2019Assortment} propose a multicut outer-approximation algorithm that works in a cutting plane fashion  by solving an MILP at every iteration. This algorithm generates cuts for groups of demand points instead of one cut of every demand point as in \cite{Ljubic2017} or one cut for all demand points as in the classical outer-approximation scheme \citep{OA_Duran1986outer,OA_Bonami2008algorithmic}. The branch-and-cut proposed by  \citep{Ljubic2018outer} and multicut outer approximation are considered as state-of-the art approaches for the MCP under the MNL model. 
Note that in the context of the MCP, MNL and MMNL instances have similar structures. Thus, all the methods developed for the MNL model can be generally applied to MMNL problem instances. There are also a couple of studies investigating the MCP under the MMNL model  \citep{Haase2009,Haase2013}. These studies make use of MILP formulations, which is generally outperformed by the  branch-and-cut approached \citep{Ljubic2018outer}.

\noindent
\textbf{Paper outline:}
The rest of paper is structured as follow.
Section \ref{sec:MCP-and-GEV} briefly presents the GEV family focusing on the some essential properties of the CPGF, and the MCP under  the GEV family. In Section \ref{sec:MCP-uder-GEV}, we investigate the monotonicity and submodularity of the MCP under the GEV family, and present our local search algorithm. Section \ref{sec:Expr.} reports computational results. Finally, Section \ref{sec:conc} concludes.

\noindent
\textbf{Notation:}
Boldface characters represent matrices (or vectors), and $a_i$ denotes the $i$-th element of vector $\ba$. We use $[m]$, for any $m\in \mathbb{N}$, to denote the set $\{1,\ldots,m\}$. 

%%%%%%%%%%%%%%%%%%%%%%%%%%%%%%%
%%%%%%%%%%%%%%%%%%%%%%%%%%%%%%%

\section{Generalized Extreme Value Models and the Maximum Capture Problem}
\label{sec:MCP-and-GEV}
In this section  we introduce some basic concepts and properties of the GEV family and formulate the MCP under GEV models. 

\subsection{Generalized Extreme Value Models}
The Random Utility Maximization (RUM) framework \citep{McFa78} is the most popular approach to  model discrete choice behavior. Under the RUM principle,  the decision maker is assumed to associate an utility $u_{j}$ with each alternative/option $j$ in a given choice set $S$ that contains all possible alternatives.
The additive RUM \citep{McFa78,FosgBier09} assumes that each random utility is  a sum of two parts $u_{j} = v_{j}+\epsilon_{j}$, where the term $v_{j}$ is deterministic and can include values  representing characteristics of the alternative and/or the decision maker, and the random term $\epsilon_{j}$  is unknown to the analyst. There are several assumptions that have been made on the randoms terms, which leads to  different types of discrete choice models in the literature, e.g., the MNL or nested logit models  \citep{McFa78,Trai03}.  
The deterministic terms $v_{j}$ often have a linear structure, i.e.,  $v_{j} = \beta^\transpose \alpha_{j}$, where $^\transpose$ is the transpose operator and $\beta$ is a vector of parameters to be estimated from historical data of how people make decisions, 
and $\alpha_{j}$ is a vector of attributes of alternative $j$.
The RUM principle then assumes that a decision is made by maximizing the random utilities, and the probability that an alternative $j$ is selected can be computed as $P(u_j\geq u_k,\;\forall k\in S)$. 

The GEV family covers most of the existing discrete choice models in the literature. This family of model is fully flexible, in the sense that it allows to construct various discrete choice models that are consistent with the RUM principle \citep{McFa81}.
Assume that the choice set contains $m$ alternative indexed as $\{1,\ldots,m\}$ and let $\bU = \{v_1,\ldots,v_m\}$ be the vector of utilities . 
A GEV model can be determined by a choice probability generating function (CPGF) $G(\bY)$ \citep{McFa81,Fosgerau2013}, where $\bY$ is a vector of size $m$ with entries $Y_j = e^{v_j}$. Given $j_1,\ldots, j_k \in  [m]$, let $\partial G_{j_1...j_k}$,  be the mixed partial derivatives of $G$ with respect to $Y_{j_1},\ldots,Y_{j_k}$.  
It is well-known that 
the CPGF  $G(\cdot)$ and the mixed partial derivatives have the
the following properties \citep{Mcfadden1978modeling}.
\begin{remark}\label{propert:GEV-CPGF}
{\it A CPGF $G(\bY)$ of a GEV model, has the following properties.
\begin{itemize}
 \item[(i)] $G(\bY) \geq 0,\ \forall \bY\in \bbR^m_+$,
 \item[(ii)] $G(\bY)$ is homogeneous of degree one, i.e., $G(\lambda \bY) = \lambda G(\bY)$
 \item[(iii)] $G(\bY)\rightarrow \infty$ if $Y_j\rightarrow \infty$
 \item[(iv)]  Given $j_1,\ldots,j_k \in [m]$ distinct from each other,
  $\partial G_{j_1,\ldots,j_k}(\bY)>0$ 
 if $k$ is odd, and $\leq$ if $k$ is even
 \item[(v)] $G(\bY) = \sum_{j\in [m]} Y_j\partial G_j(\bY)$
 \item[(vi)] $\sum_{k\in [m]} Y_k\partial G_{jk} (\bY) = 0$, $\forall j\in [m]$.
\end{itemize}
%where $\partial G_i(\bY) = \partial G(\bY)/\partial Y_i$
}
\end{remark}  
Here we note that $(i)-(iv)$ are basic properties of a GEV generating function \citep{McFa81}, and Properties $(v)$ and $(vi)$ are direct results from the homogeneity property.
We will make use of these properties throughout the rest of the paper to explore the properties of the objective function of the MCP, and derive solution algorithms for the MCP.
Under a GEV model specified by a CPGF $G(\bY)$, the choice probability of an alternative $j\in [m]$, conditional on $\bY$ and $G$, is given by 
\[
P(j|\bY,G)  = \frac{Y_i \partial G_i(\bY)}{G(\bY)}.
\]
The GEV framework allows for correlated utilities and one can build different CPGF to model different correlation patterns among random utilities. One can build a GEV model from a network of correlation structure, which provides a very flexible way to construct choice models that are able to capture complex relationships between alternatives \citep{MaiFreFosBas15_DynMEV,DalyBier06}. In the following, we show some specific instances of the GEV family that are already popular in the demand modeling and operations research literature.  

\subsubsection*{The MNL model.}
 The MNL is one of the most widely-used models in the literature.
%Its popularity may be because of its simple structure. 
 This model results from the   assumption that the random terms $\epsilon_j$, $j\in [m]$, are independent and
identically distributed (i.i.d.) and follow the standard Gumbel distribution. The CPGF function has a simple form as $G(\bY) = \sum_{j\in [m]}Y_j$ and the choice probabilities have the fractional form below
\begin{align}
    P(j|\bY,G) = \frac{Y_j}{\sum_{j\in [m]} Y_j} = \frac{e^{v_j}}{\sum_{j\in [m]} e^{v_j}}.
\end{align}
It is well-known that the MNL model exhibits from the IIA property, which means that the choice probability of an alternative will not be affected by the attributes or the state of the other alternatives. However, in some situations, alternatives share unobserved attributes (i.e, random terms are correlated) and the IIA property does not hold.

\subsubsection*{The nested logit model.}
The nested  logit model \citep{BenA73} is one of the first attempts to relax the IIA property from the  MNL model. 
In this GEV model, the choice set is partitioned into $L$ nests, which are disjoint subsets of alternatives. Let denote by $n_1,\ldots,n_L$ the $L$ nests.
The corresponding  CPGF can be written as 
\[
G(\bY) = \sum_{l\in L} \left(\sum_{j\in n_l} Y_j^{\mu_l} \right)^{1/\mu_l}
\]
where  $\mu_l\geq 1$, $l\in[L]$, are the parameters of the nested model. This model is based on the observation that, in many situations, some similar or closely related  alternatives can be grouped into smaller subsets. It is easy to see that the function $G$ above satisfies the six properties above and the choice probabilities can be computed as 
\[
P(j|\bY,G) = \frac{\left(\sum_{j'\in n_l} Y_{j'}^{\mu_l}\right)^{1/\mu_l}}{\sum_{l\in [L]}\left(\sum_{j'\in n_l} Y_{j'}^{\mu_l}\right)^{1/\mu_l}} \frac{Y_j^{\mu_l}}{\sum_{j'\in n_l} Y_{j'}^{\mu_l}},\; \forall l\in [L], j\in n_l.  
\]
The cross-nested logit model  \citep{BenABier99a} is an extension of the nested logit that allows the nests to share common alternatives. This model is known to be fully flexible, as it can approximate arbitrarily close any RUM models \citep{Fosgerau2013}. The network GEV model proposed in \cite{DalyBier06} further generalizes the cross-nested model by proving a way to construct a GEV CPGF based on any rooted network of correlation structure.    

Beside the GEV family, the \textbf{MMNL model} \citep{McFaTrai00} is also popular due to its flexibility in capturing utility correlation. In the MMNL model, the model parameters (and the utilities $v_j$) are assumed to be random, and the choice probabilities can be obtained by taking the expectations over random coefficients. Let $\bY^1,\ldots,\bY^K$ be $K$ realizations sampled from  the distribution of the random parameters, the choice probabilities can be approximated as
\[
P(j|\bY^1,\ldots, Y^K,G) =  \frac{1}{K}\sum_{k=1}^K \frac{Y_j^k}{\sum_{t\in [m]} Y^k_t}.
\]
The MMNL model is highly preferred in practice due to its flexibility in modeling people demand. However, the estimation and application of this model in decision-making are well-known to be expensive and complicated, due to the fact that it  requires simulation to approximate the choice probabilities.

\subsection{The Maximum Capture Problem}

We are interested in the situation that  a ``\textit{newcomer}'' firm  wants to locate new facilities in a competitive market, i.e., there are already existing facilities from competitors that can serve  customers. The firm may want to maximize the expected market share achieved by attracting the customers to new facilities. To capture the customers' demand, we suppose that a customer selects a facility according to a RUM model. In this context, each customer associate each facility  with a random utility and we  assume that the  customer will choose a facility by maximizing his/her utilities. Accordingly, the firm aims at selecting a set of locations to locate new facilities to maximize the expected number of customers. 
In the following, we describe in detail the MCP under GEV models. 

We denote by $\cV = [m]$ the set of possible locations.
 Let $I$ be the set of geographical zones where customers are located and $q_{i}$ is the number of customers in zone $i\in I$  and for customers at zone $i$, let $v_{ij}$ be the corresponding deterministic utility of location $j\in [m]$. These utility values can be inferred by estimating the RUM model using historical data.  The set $I$ can be viewed as a set of customer types, e.g., customers that belong to different categories specified by, for instance, age or income.   
A GEV model for customers located at zone $i\in I$ can be represented by a choice probability generating function (CPGF) $G^i(\bY^i)$, where $\bY^i$ is a vector of size $m$ with entries $Y^i_j = e^{v_{ij}}$. 

%Given $j_1,\ldots, j_k \in  [m]$, let $\partial G^i_{j_1...j_k}$,  be the mixed partial derivatives of $G$ with respect to $Y^i_{j_1},\ldots,Y^i_{j_k}$.  
%In the following, we describe the problem with three different choice models in detail, i.e, the MNL. the MMNL, and the nested logit model. 
Under a GEV model specified by a set of CPGF $G^i(\bY^i)$, $i\in I$, taking into consideration the competitors, the choice probability  of a location $j\in [m]$  is given as 
\[
P(j|\bY^i,G^i) = \frac{Y_j \partial G^i_j(\bY^i)}{1+ G^i(\bY^i)}.
\]
Here, without loss of generality,  we assume that the  total utility of the competitor is 1 for the sake of simplicity. 
As if it is not the case, then we can always scale the utilities $\bY^i$ to get utilities of $1$ for the competitors. More specifically, it is possible due to fact that, for any $\alpha>0$, 
\begin{align}
\frac{Y_j \partial G^i_j(\bY^i)}{\alpha + G^i(\bY^i)} &\stackrel{(a)}{=} \frac{\frac{Y_j}{\alpha}  \partial G^i_j(\bY^i)}{1 + G^i({\bY^i}/{\alpha})} \stackrel{(b)}{=} \frac{\frac{Y_j}{\alpha}  \partial G^i_j(\bY^i/\alpha)}{1 + G^i({\bY^i}/{\alpha})} \nonumber 
\end{align}
where $(a)$ is due the homogeneity of $G^i(\cdot)$ ($(ii)$ of Remark \ref{propert:GEV-CPGF}) and $(b)$ is obtained by taking derivatives of the both sides of the equation $G^i(\alpha \bY^i) = \alpha G^i(\bY^i)$ w.r.t. $Y^i_j$
\[
\alpha \partial G^i_j(\alpha \by^i) = \alpha \partial G^i_j(\bY^i),\text{ or } \partial G^i_j(\alpha \by^i) = \partial G^i_j(\bY^i),\;\text{ for any }  \alpha >0. 
\]
We are interested in the fact that the facilitates are located at a subset of locations $S\subset [m]$. 
Hence, the conditional choice probability can be written as
\[
P(j|\bY^i,G^i, S) = \frac{Y^i_j \partial G^i_j(\bY^i|S)}{1+ G^i(\bY^i|S)},\; \forall j\in S,
\]
where the conditional CPGF $G^i(\bY^i|S)$ can be computed as $G^i(\bY^i|S) = G^i(\widetilde{\bY}^i)$, where $\widetilde{\bY}^i$ is a vector of size $m$ with entries $\widetilde{Y}^i_j = {Y}^i_j$ if $j\in S$ and $\widetilde{Y}^i_j =0$ otherwise. This can be interpreted as if a location $j$ is not in $S$, then its corresponding utility becomes  very small, i.e., $v_{ij} = -\infty$, then $Y^i_j = e^{v_{ij}} = 0$.  
The maximum capture problem under a GEV model specified by CPGFs $G^i(\bY^i)$, $i\in I$, can be stated as 
\begin{equation}\label{prob:MCP-1}
 \max_{S \in \cS}\left\{f^{\GEV}(S) = \sum_{i\in I}q_i\sum_{j\in S} P(j|\bY^i,G^i, S) \right\},
\end{equation}
where $\cS$ is the set of feasible solutions. Under a cardinality constraint $|S| \leq C$, $\cS$ can be defined as $\cS = \{\cS\subset[m]|\; |S|\leq C\}$, for a given constant $C$ such that $1\leq C\leq m$. 
Note that the objective function can be further simplified as
\begin{align}
    f^{\GEV}(S) &=  \sum_{i\in I} q_i \frac{\sum_{j\in S} Y^i_j \partial G^i_j(\bY^i|S)}{ 1+ G^i(\bY^i|S)}\nonumber \\
    & \stackrel{(a)}{=}  \sum_{i\in I} q_i - \sum_{i\in I}  \frac{q_i}{ 1+ G^i(\bY^i|S)},\nonumber
\end{align}
where $(a)$ is due to Property $(v)$ in Remark \ref{propert:GEV-CPGF}. 

If the choice model is MNL, the objective function becomes 
\[
f^{\GEV}(S) = \sum_{i\in I} q_i - \sum_{i\in I}  \frac{q_i}{ 1+ \sum_{j\in m} Y^i_j},
\]
and from previous studies, we know that $f^{\GEV}(S)$ is submodular and it binary representation is concave \citep{Benati2002}. Thus, an approach based on sub-gradient and submodular cuts can be used \citep{Mai2019Assortment,Ljubic2018outer} to efficiently solve the problem. 
On the other hand, formulations based on GEV models would be much more complicated. For example, under the nested logit model, we can write the objective function as
\[
f^{\GEV}(S) =  \sum_{i\in I} q_i - \sum_{i\in I}  \frac{q_i}{ 1+ \sum_{l\in [L]} \left(\sum_{j\in n_l\cap S} (Y^i_j)^{\mu_l}\right)^{1/\mu_l} }.
\]
 Under a general case, e.g., the network GEV model \cite{DalyBier06, MaiFreFosBas15_DynMEV}, it is even not possible to write the objective function in a closed form.

It is important to note that, if we look at the objective function under a MMNL model
\[
\begin{aligned}
f^{\MMNL}(S) &= \frac{1}{K}\sum_{k\in[K]} \left( \sum_{i\in I} q_i - \sum_{i\in I}  \frac{q_i}{ 1+ \sum_{j\in m} Y^{i,k}_j}\right) \\
&=   \sum_{k\in [K],i\in I} \frac{q_i}{K} - \sum_{i\in I,k\in [K]}  \frac{q_i/K}{ 1+ \sum_{j\in m} Y^{i,k}_j},
\end{aligned}
\]
where $\bY^{i,k}$, $k=1,...,K$, are $K$ realization of the random utility vector $\bY^i$. So, this objective function can be viewed as one from the MNL-based MCP problem with  $K\times |I|$ customer zones, in which there are $q_i/K$ customers in zone $(i,k)$-th. So, all the results established for the MNL (and GEV in general) problem
can also be used to solve the MMNL problem.

\section{Maximum Capture Problem under GEV Models}
\label{sec:MCP-uder-GEV}
 In this section we explore the MCP under GEV models. In particular, by leveraging the properties of the GEV CPGFs shown above, 
 we show that the objective function in the context is monotonic and submodular.
 %, which provides a performance guarantee for a simple greedy heuristic, and  
 This generalizes some well-known results established for MNL-based problems in previous studies \citep{Benati2002}. We also design a location search procedure  to efficiently solve the problem.  
%%%%%%%%%%%%%%%%%%%%%%%%%%%%

\iffalse
\subsection{Problem Formulation}
We denote by $\cV = [m]$ the set of possible locations. A GEV model for customers located at zone $i\in I$ can be represented by a choice probability generating function (CPGF) $G^i(\bY^i)$, where $\bY^i$ is a vector of size $m$ with entries $Y^i_j = V_{ij}$. Given $j_1,\ldots, j_k \in  [m]$, let $\partial G^i_{j_1...j_k}$,  be the mixed partial derivatives of $G$ with respect to $Y^i_{j_1},\ldots,Y^i_{j_k}$.  
It is well-known that 
the CPGF  $G^i(\cdot)$ and the mixed partial derivatives have the
the following properties \citep{Mcfadden1978modeling}.
\begin{remark}\label{propert:GEV-CPGF}
{\it A CPGF $G^i(\bY^i)$, $i\in I$, has the following properties.
\begin{itemize}
 \item[(i)] $G^i(\bY^i) \geq 0,\ \forall \bY^i\in \bbR^m_+$,
 \item[(ii)] $G^i(\bY^i)$ is homogeneous of degree one, i.e., $G^i(\lambda \bY^i) = \lambda G^i(\bY^i)$
 \item[(iii)] $G(\bY^i)\rightarrow \infty$ if $Y^i_j\rightarrow \infty$
 \item[(iv)]  Given $j_1,\ldots,j_k \in [m]$ distinct from each other,
  $\partial G^i_{j_1,\ldots,j_k}(\bY^i)>0$ 
 if $k$ is odd, and $\leq$ if $k$ is even
 \item[(v)] $G^i(\bY^i) = \sum_{j\in [m]} Y^i_j\partial G^i_j(\bY^i)$
 \item[(vi)] $\sum_{k\in [m]} Y^i_k\partial G_{jk} (\bY^i) = 0$, $\forall j\in [m]$.
\end{itemize}
%where $\partial G_i(\bY) = \partial G(\bY)/\partial Y_i$
}
\end{remark} 
\fi

%%%%%%%%%%%%%%%%%%%%%%%%%%%%%%%%%
%%%%%%%%%%%%%%%%%%%%%%%%%%%%%%%%%
\subsection{Monotonicity and Submodularity}

We show two key results of the paper, which demonstrates that the objective function under any GEV model is monotonic and submodular, thus providing a performance guarantee  for a greedy heuristic procedure.

To prove the results, we first formulate the MCP as a binary program. That is, given a subset $S\subset [m]$, let $\bx^S$ be a binary vector of size $m$ with entries $\bx^S_j = 1$ if  $j\in S$ and $\bx^S_j = 0$ otherwise. We see that the conditional CPGF  can be written  as
\[
G^i(\bY^i|S) = G^i(\bx^S \circ \bY^i), 
\]
where $\circ$ is the element-by-element operator and $\bx \circ \bY^i$ is vector of size $m$ with entries $x^S_j Y^i_j$, $j=1,\ldots,m$. We now can formulate  \eqref{prob:MCP-1} as
\begin{equation}\label{prob:MCP-2}
 \max_{\bx \in X }\left\{f^{\GEV}(\bx) = \sum_{i\in I}q_i   - \sum_{i\in I} \frac{q_i}{1+ G^i(\bx \circ \bY^i)} \right\},
\end{equation}
where $X  = \{\bx^S \in\{0,1\}^m|\; \forall S\in \cS  \}$, i.e.,  the feasible set of binary solutions that corresponds to all the subsets in $\cS$. It is worth noting that if the choice model is MNL (or MMNL), then the objective $f^{\GEV}(\bx)$ is concave in $\bx$ and the problem can be handled efficiently by an outer-approximation method  \citep{Bonami2011_BB_MIP,Mai2019Assortment,Ljubic2018outer}. It is however not the case under an arbitrary  GEV model.

The following proposition tells us that the objective function is monotonic, which implies that adding more facilities always yields better objective values. 
\begin{proposition}[Monotonicity]
\label{prop:adding-new-facilities}
Adding more facilities always yields better objective values,
i.e., $f^{\GEV}(S\cup \{i\}) > f^{\GEV}(S)$ for any $i\notin S$.
\end{proposition}
\proof{Proof:}
To prove the claim, let $\bx \in \{0,1\}^m$ be the binary vector  representing set $S$. For any $j\in [m]$ such that  $x_j = 0$, we need to prove $f^{\GEV}(\bx + \bbe^j) > f^{\GEV}(\bx)$, where $\bbe^j$ is a vector of size $m$ with zero entries except the $j$-th element that is equal to 1. To prove this, let us consider $G^i(\bx \circ \bY^i)$. Taking the derivative of $G^i(\bx \circ \bY^i)$ w.r.t. an $x_j$, $j\in [m]$, we have
\begin{align}
    \frac{\partial G^i(\bx \circ \bY^i)}{ \partial x_j} = Y^i_j \partial G^i_j (\bx \circ \bY^i) \stackrel{(b)}{>}0,
\end{align}
where $(b)$ is due to Property $(iv)$ of Remark \ref{propert:GEV-CPGF}. This implies that $G^i(\bx \circ \bY^i)$ is (strictly) monotonic increasing in any $x_j$, $j\in [m]$. Thus   
\begin{equation}
\label{eq:prop1-eq1}
G^i((\bx+ \bbe^j) \circ \bY^i) > G^i(\bx \circ \bY^i).
\end{equation}
Together with the definition of $f^{\GEV}(\bx)$ in \eqref{prob:MCP-2}, we have $f^{\GEV}(\bx + \bbe^j) > f^{\GEV}(\bx)$ as desired. 
\endproof
Since $f^{\GEV}(S)$ is monotonic,   if we consider a cardinality constraint $|S|\leq C$ for a scalar $C \in \{1,\ldots,m\}$, then an optimal solution $S^*$ always achieves the maximum cardinality, i.e., $|S^*| = C$. Thus, we can replace the cardinality constraint by an equality one, i.e., $|S| = C$. Note that a similar claim has been validated for the MNL-based problems in prior work  \citep{Mai2019Assortment}. 

The submodularity is well-known for the objective function under the MNL model \citep{Benati2002}. The theorem below shows that it is also the case under any models in the GEV family. 
\begin{theorem}[Submodularity]
\label{th:th1}
$f^{\GEV}(S)$  is submodular.
\end{theorem}
\proof{Proof:}
To prove the submodularity, we will show that for any set $A\subset B \subset [m]$  and for any $j \in [m]\backslash  B$ we have
\begin{equation}\label{eq:proof-th1-eq1}
f^{\GEV}(A \cup \{j\})  - f^{\GEV}(A) \geq f^{\GEV}(B \cup \{j\}) -  f^{\GEV}(B)  
\end{equation}
Let use denote each component of $f^{\GEV}(S)$ as
\[
g^i(S) = \frac{q_i}{ 1+  G^i(\bY^i|S)}, \forall i\in I 
\]
then \eqref{eq:proof-th1-eq1} can be validated if we can prove
\[
g^i(A \cup \{j\})  - g^i(A) \leq g^i(B \cup \{j\}) -  g^i(B),
\]
or equivalently,
\begin{align}
    &\frac{1}{1+ G^i((\bx^A + \bbe^j) \circ \bY^i)} - \frac{1}{1+ G^i(\bx^A \circ \bY^i)} \leq    \frac{1}{1+ G^i((\bx^A + \bbe^j) \circ \bY^i)} - \frac{1}{1+ G^i(\bx^A \circ \bY^i)} \nonumber \\
    &\Leftrightarrow  \frac{G^i((\bx^A + \bbe^j) \circ \bY^i) - G^i((\bx^A) \circ \bY^i)}{[1+ G^i((\bx^A + \bbe^j) \circ \bY^i)][1+ G^i(\bx^A \circ \bY^i)]} \geq \frac{G^i((\bx^B + \bbe^j) \circ \bY^i) - G^i(\bx^B \circ \bY^i)}{[1+ G^i((\bx^B + \bbe^j) \circ \bY^i)][1+ G^i(\bx^B \circ \bY^i)]}  \label{eq:proof-th1-eq2}
\end{align}
Now, for ease of notation, for any $\bx\in \{0,1\}^m$ and $j\in [m]$ such that $x_j = 0$,  let 
\[
\phi(\bx) =  G^i((\bx + \bbe^j) \circ \bY^i) - G^i(\bx \circ \bY^i)
\]
For any $k\in [m]$ such that $x_k = 0$ and $k\neq j$ we take the partial derivative of $\phi(\bx)$ w.r.t. $x_k$ and get
\begin{equation}\label{eq:proof-th1-eq3}
\frac{\partial \phi(\bx)}{ \partial x_k} = Y^i_k \left( \partial G^i_k((\bx + \bbe^j) \circ \bY^i) - \partial G^i_k(\bx \circ \bY^i)\right)
\end{equation}
Now, let define another function $\rho(\bx) = \partial G^i_k(\bx \circ \bY^i)$. Taking the partial derivative of $\rho(\bx)$ w.r.t. $x_j$ we get
\[
\frac{\partial \rho(\bx)}{\partial x_j} = Y^i_j \partial G^i_{kj} (\bx \circ \bY^i),
\]
and since $\partial G^i_{kj} (\bx \circ \bY^i) \leq 0$ (Property $(iv)$ of Remark \ref{propert:GEV-CPGF}), we have 
${\partial \rho(\bx)}/{\partial x_j}\leq 0$. Thus, $\rho(\bx)$ is monotonic decreasing in $x_j$, which implies
\begin{equation}
\label{eq:proof-th1-eq4}
    \partial G^i_k((\bx + \bbe^j) \circ \bY^i) - \partial G^i_k(\bx \circ \bY^i) \leq 0.
\end{equation}
Combine \eqref{eq:proof-th1-eq3} and \eqref{eq:proof-th1-eq4} we have $\partial \phi(\bx)/\partial x_k \leq 0$. Thus, $\phi(\bx)$ is monotonic decreasing in $x_k$, leading to the inequality
\[
 G^i((\bx + \bbe^j) \circ \bY^i) - G^i(\bx \circ \bY^i) \geq  G^i((\bx + \bbe^j + \bbe^k) \circ \bY^i) - G^i((\bx+\bbe^k) \circ \bY^i).
\]
Consequently, we have 
\begin{equation}
\label{eq:proof-th1-eq5}
    G^i((\bx^A + \bbe^j) \circ \bY^i) - G^i((\bx^A) \circ \bY^i) \geq G^i((\bx^B + \bbe^j) \circ \bY^i) - G^i((\bx^B) \circ \bY^i),
\end{equation}
for any $A\subset B\subset [m]$ and $j\notin B$. Moreover, using \eqref{eq:prop1-eq1} from the proof of Proposition \ref{prop:adding-new-facilities}, since  $A\subset B$,  we have
\begin{align}
G^i(\bx^A \circ \bY^i) &\leq G^i(\bx^B \circ \bY^i) \nonumber \\
    G^i((\bx^A + \bbe^j)\circ \bY^i) &\leq G^i((\bx^B+ \bbe^j)\circ \bY^i). \nonumber 
\end{align}
Thus, 
\begin{equation}
\label{eq:proof-th1-eq6}
[1+ G^i((\bx^A + \bbe^j) \circ \bY^i)][1+ G^i(\bx^A \circ \bY^i)] \leq [1+ G^i((\bx^B + \bbe^j) \circ \bY^i)][1+ G^i(\bx^B \circ \bY^i)]
\end{equation}
Combine \eqref{eq:proof-th1-eq5} and \eqref{eq:proof-th1-eq6} we obtain \eqref{eq:proof-th1-eq2} and then \eqref{eq:proof-th1-eq1} as desired. 
\endproof
Theorem \ref{th:th1}  generalizes the submodularity of the MNL-based MCP problem shown in previous studies \citep{Benati2002}. Together with the fact that $f^{\GEV}(S)$ is monotonic (Proposition \ref{prop:adding-new-facilities}), we know that a simple greedy local search algorithm will guarantee an $(1-1/e)$ approximation solution, i.e.,  a greedy will return a solution $\overline{S}$ such that $f^{\GEV}(\overline{S}) \geq (1-1/e)\max_{S\in \cS} f^{\GEV}(S)$ \citep{Nemhauser1978analysis}.    
\begin{corollary}[Performance guarantee for a greedy heuristic]
Under a cardinality constraint, a greedy heuristic algorithm can guarantee an $(1-1/e)$ approximation solution.
\end{corollary}

%%%%%%%%%%%%%%%%%%%%%%%%%%%%%%
%%%%%%%%%%%%%%%%%%%%%%%%%%%%%%
\subsection{Gradient-based Local Search}

Due to the submodularity property, a greedy heuristic can guarantee an $(1-1/e)$ approximation solution. In this section, we design a new local search procedure to further improve this greedy solution. Our approach is motivated by the fact  that the objective function $f^{\GEV}(\bx)$ is differentiable, suggesting that we could use gradient information to direct the search. 
The general idea to design an iterative procedure, where at each step we build a model function (linear or quadratic) to approximate the objective function using gradient and/or Hessian information. We then maximize the  model function to find a new iterate. A key component of our approach is that the model function can be only  an adequate representation of the objective function in a local neighbourhood of the current solution. Thus, we  only maximize the model function within a restricted region.
%, and the size of region can be adjusted according to the accuracy of the model function.
This approach is inspired by  the trust-region method  widely used in continuous optimization \citep{Conn2000trust}.

To start our exposition, let us define a model function based on Taylor series built around a solution candidate $\overline{\bx}$
\[
f^{\GEV}(\bx) \approx  f^{\GEV}(\overline{\bx}) + \nabla f^{\GEV}(\overline{\bx})^{\transpose} (\bx - \overline{\bx}) + \frac{1}{2}(\bx - \overline{\bx})^\transpose \bB (\bx - \overline{\bx}),  
\]
where $\bB$ is the Hessian matrix or an approximation of it at $\overline{\bx}$. In our context, the Hessian can be computed easily, but  maximizing the model function will involve solving a binary quadratic maximization problem, which is expensive. Thus, we set $\bB_k = 0$. In other words, we use a linear model function to approximate $f^{\GEV}(\bx)$. 

At each iteration, we need to solve the following sub-problem
\begin{align}
	\underset{\bx}{\text{max}}\qquad &  \nabla f^{\GEV}(\overline{\bx})^{\transpose} \bx  & \label{prob:sub-P-x}\tag{P1} \\
	 \text{subject to} \qquad & \sum_{j\in [m]} x_j  = C &  \label{cstr:sub-P-x-1}\\
	  &  \sum_{j\in [m]} |x_j - \overline{x}_j|\leq \Delta  &  \label{cstr:sub-P-x-2}\\
	  &  \bx \in \{0,1\}^m &  \nonumber
\end{align}
where \eqref{cstr:sub-P-x-1} is the cardinality constraint, and \eqref{cstr:sub-P-x-2} is to ensure that the new solution candidate is within a region of size $\Delta$ around $\overline{\bx}$. Note that \eqref{cstr:sub-P-x-2} can  be linearized as
\[
\sum_{j\in[m], \overline{x}_j = 1} (1-x_j)  + \sum_{j\in[m], \overline{x}_j = 0} x_j \leq \Delta ,
\]
so as \eqref{prob:sub-P-x} becomes a integer linear program, which can be handled by an MILP solver. In the following, we will look closely to \eqref{prob:sub-P-x} and show that it can be solved to optimality in polynomial time.  

%%%%%%%%%%%%%%%%%%%%%%%%%%%%%%
%%%%%%%%%%%%%%%%%%%%%%%%%%%%%%
\noindent\textbf{{Solving Subproblems:}}
We further look into the subproblem of the gradient-based local search  \eqref{prob:sub-P-x} to design an efficient algorithm to solve it.
To facilitate our exposition, we first note that the constraint $\sum_{j\in [m]}|\overline{x}_j - x_j| \leq \Delta$ implies that there are at most $\Delta/2$ locations that either appears in $S$ or in $\overline{S}$, but not in both, where $S,\overline{S}$ are the subsets representing $\bx$ and $\overline{\bx}$, respectively. 
For this reason, $\Delta$ should be integer and even, and  the constraint $\sum_{j\in [m]}|\overline{x}_j - x_j| \leq \Delta$ is equivalent to $|S \bigtriangleup \overline{S}| \leq \Delta$, where $\bigtriangleup$ is the symmetric difference operator, i.e., $S \bigtriangleup \overline{S} = (S\backslash \overline{S}) \cup (\overline{S}\cup S)$. We therefore can rewrite \eqref{prob:sub-P-x} as 
\begin{align}
	\underset{S \subset [m]}{\text{max}}\qquad &  \sum_{j\in S} d_j  & \label{prob:sub-P-S}\tag{P2} \\
	 \text{subject to} \qquad & |S|  = C&  \label{cstr:sub-P-S-1}\\
	  &   |S \bigtriangleup \overline{S}| \leq \Delta, &  \label{cstr:sub-P-S-2}
\end{align} 
where $\overline{S}\subset [m]$ is the subset that corresponds to the binary vector $\overline{\bx}$ and $d_j = \nabla f^{\GEV}(\overline{x})_j$, $j\in [m]$. Under the cardinality constraint $|S| = C$, we  see that $|S\backslash \overline{S}| = |\overline{S} \backslash {S}| =  |S \bigtriangleup \overline{S}| / 2$. 
The following proposition shows that $d_j$ are non-negative, for all $j\in [m]$. 
\begin{proposition}
All the coefficients of the objective function of \eqref{prob:sub-P-S} are non-negative.
\end{proposition}
\proof{Proof:}
We prove the claim by showing that, for  any $\bx \in [0,1]^m$, $\nabla_\bx f^{\GEV}(\bx) \geq 0$.
Given $j\in [m]$,
 taking the derivative of $f^{\GEV}(\bx)$ w.r.t. $x_j$ we have
\begin{align}
\frac{\partial f^{\GEV}(\bx) }{\partial x_j} &= \sum_{i\in I} \frac{\partial G^i(\bx \circ \bY^i)}{\partial x_j} \frac{q_i}{(1+G^i(\bx \circ \bY^i))^2} \nonumber \\
&= \sum_{i\in I} \frac{ q_iY^i_j  \partial G^i_j(\bx \circ \bY^i)}{(1+G^i(\bx \circ \bY^i))^2} \geq 0 \label{eq:prop2-eq1}
\end{align}
where \eqref{eq:prop2-eq1} is due to the fact that   $\partial G^i_j(\bx \circ \bY^i)>0$ (Property $(iv)$ of Remark \ref{propert:GEV-CPGF}). We obtain the desired inequality.
\endproof

In Algorithm \ref{algo:sub-problem} we describe our main steps to solve \eqref{prob:sub-P-S}. In Step 1, we find $\Delta/2$ smallest coefficients $d_j$ in $\overline{S}$ and $\Delta/2$ largest coefficients $d_j$ in $[m]\backslash \overline{S}$. This is motivated by the fact we only seek subsets generated by exchanging at most $\Delta/2$ elements  in $\overline{S}$ with some outside $\overline{S}$. Thus, to achieve best objective values, we should exchange elements of lowest coefficients in $\overline{S}$ with those of highest coefficients in $[m]\backslash \overline{S}$. In the second step, $\gamma(t)$ represents the best gain obtained by exchanging $t$ elements, and in the third step we just select the best $\gamma(t)$ to get the best solution. Proposition \ref{prop:convergence-sub-P-S} below shows that Algorithm \ref{algo:sub-problem} will efficiently return an optimal solution to \eqref{prob:sub-P-S}.
\begin{algorithm}[htb]
    \caption{Solving sub-problems} \label{algo:sub-problem}
    \comments{Step 1: Take smallest coefficients in $\overline{S}$ and largest coefficients in $[m]\backslash \overline{S}$}\\
    Choose $\sigma^1_1,\ldots,\sigma^1_{\Delta/2} \in \overline{S}$ and $\sigma^2_1,\ldots,\sigma^2_{\Delta/2}\in [m]\backslash \overline{S}$ such that
    \[
    \begin{aligned}
 d_{\sigma^1_1}\leq d_{\sigma^1_2} \leq \ldots \leq d_{\sigma^1_{\Delta/2}} \leq \min_{j\in \overline{S}\backslash \{\sigma^1_1,...,\sigma^1_{\Delta/2}\}} d_j  \\
 d_{\sigma^2_1}\geq d_{\sigma^2_2} \geq \ldots \geq d_{\sigma^2_{\Delta/2}} \geq \max_{j\in [m]\backslash \overline{S} \backslash \{\sigma^2_1,...,\sigma^2_{\Delta/2}\}} d_j  \\
    \end{aligned}
    \]
    \comments{Step 2: Select the best set for each local region size $|S \bigtriangleup \overline{S}|  = 2t,\; \text{ for } t=1,\ldots,\Delta/2$}\\
    \For{$t=1,\ldots,\Delta/2$}
    {
   $$\gamma(t) = \sum_{h=1}^t \left(d_{\sigma^2_h} - d_{\sigma^1_h} \right)$$
    }
    Select $t^* = \text{argmax}_{t = 1,\ldots,\Delta/2 } \gamma(t)$ \\
    \comments{Step 3: Return the best solution within the local region $|S \bigtriangleup \overline{S}| \leq \Delta$}\\
    Return
    \[
   S^* \leftarrow  \overline{S} \cup \{\sigma^2_{1},\ldots, \sigma^2_{t^*}\} \backslash \{\sigma^1_{1},\ldots, \sigma^1_{t^*}\} 
    \]
\end{algorithm}

\begin{proposition}
\label{prop:convergence-sub-P-S}
Algorithm \ref{algo:sub-problem} returns an optimal solution to \eqref{prob:sub-P-x} with complexity $\mathcal{O}(m\Delta/2)$.
\end{proposition}
\proof{Proof:}
To prove the convergence, we let $S$ be a feasible solution of \eqref{prob:sub-P-S}, i.e., $|S| = C$  and $|S \bigtriangleup \overline{S}| \leq \Delta$. We will prove that $\sum_{j\in S} d_j \leq \sum_{j\in S^*}d_j$, where $S^*$ is the solution returned. from Algorithm \ref{algo:sub-problem}. 
Let $t = |S \bigtriangleup \overline{S}|/2$, then we know that $S$ can be obtained by exchanging $t$ elements between $S$ and $\overline{S}$. Let $\pi^1_1,\ldots,\pi^1_t$ be the indexes of $t$ elements  in $\overline{S}$ that are exchanged with  $t$ elements in $t$, indexed as $\pi^2_1,\ldots,\pi^2_t$. We have
\begin{align}
    \sum_{j\in S} d_j &= \sum_{j\in \overline{S}}d_j - \sum_{h=1}^t d_{\pi^1_h} + \sum_{h=1}^t d_{\pi^2_h} \nonumber \\
    &\stackrel{(c)}{\leq}  \sum_{j\in \overline{S}}d_j - \sum_{h=1}^t d_{\sigma^1_h} + \sum_{h=1}^t d_{\sigma^2_h} \nonumber \\
    &\stackrel{}{=}  \sum_{j\in \overline{S}}d_j + \gamma(t) \nonumber \\
    &\stackrel{(d)}{\leq}  \sum_{j\in \overline{S}}d_j + \gamma(t^*) \stackrel{}{=}  \sum_{j\in S^*}d_j, \nonumber
\end{align}
where $(c)$ is due to the way we select $\sigma^1_h$ and $\sigma^2_h$, $h=1,\ldots,\Delta/2$, and $(d)$ is due to the way $t^*$ is selected. 
This implies that $S^*$ is an optimal solution to \eqref{prob:sub-P-S} as desired. 

For the complexity, we see that Step 1 would take $\cO(\Delta/2 |\overline{S}| + \Delta/2 (m-|\overline{S}|) ) = \cO(m\Delta/2)$. Step  2 would require  $\cO(\Delta^2/4)$, which would be much smaller than $\cO(m\Delta/2)$. Adding all together, the complexity of Algorithm \ref{algo:sub-problem} is $\cO(m\Delta/2)$. 
\endproof

%%%%%%%%%%%%%%%%%%%%%%%%%%%%%%
%%%%%%%%%%%%%%%%%%%%%%%%%%%%%

\subsection{GGX Algorithm}

Our main algorithm consists of three main phases. In the first phase (warm up), we perform a greedy heuristic, which can be done by starting from the null set and adding locations one at a time, taking at each step the location that increases $f^{\GEV}(\cdot)$ the most. This phase finishes when we reach the the maximum capacity, i.e., $|S| = C$. After this phase, due to the submodularity, it is guaranteed that the obtained solution yields at least a factor $(1-1/e)$ times the optimal value. In the second phase, we iteratively solve the sub-problem \eqref{prob:sub-P-x} to  seek better solutions. This phase ends when we cannot find any better solutions. In the last phase, we further enhance the solution obtained by performing a simple greedy local search based on exchanging some locations in the current set $S$ with some others from $[m]\backslash S$. We describe the three phases in detail in Algorithm \ref{algo:local-search}. 

\begin{algorithm}[htb]
    \caption{GGX algorithm} 
    \label{algo:local-search}
    \SetKwRepeat{Do}{do}{until}
    \comments{1: \textbf{G}reedy heuristics (warm up step)}\\
    $S = \emptyset$ \\
    \For{$j=1,\ldots,C$}{
    $j^* = \text{argmax}_{j\in [m]\backslash S} f^{\GEV}(S \cup \{j\})$\\
    $S \leftarrow S \cup \{j^*\} $
    }
    \comments{2: \textbf{G}radient-based local search} \\
    $k=0; S^0 = S$ \\
    \Do{$S^k = S^{k-1}$
    }
    {
    Solve \eqref{prob:sub-P-x} based on a local region around $\bx^{S^k}$ to get a new solution candidate $\overline{S}$\\
    \If{$f^{\GEV}(\overline{S}) >f^{\GEV}(S^k)$}{$S^{k+1}\leftarrow \overline{S}$}
    \Else
    {
    $S^{k+1} = S^k$
    }
    $k \leftarrow k+1$
    }
    \comments{3: \textbf{Ex}changing phase}\\
    \Do{$S^k = S^{k-1}$
    }
    {
    $(j^*,t^*)  = \text{argmax}_{\substack{j\in S\\ t\in [m]\backslash S}} \left\{ f^{\GEV}( S^k\cup \{t\} \backslash \{j\}) \right\}$ \\
    $\overline{S} =  S^k\cup \{k^*\} \backslash \{j^*\}$\\ 
    \If{$f^{\GEV}(\overline{S}) > f^{\GEV}(S^k)$}{$S^{k+1}\leftarrow \overline{S}$}
    \Else{
    $S^{k+1} = S^k$
    }
    $k\leftarrow k+1$
    }
    Return $S^k$.
\end{algorithm}

In the context of assortment optimization under GEV models, \cite{Mai2019Assortment} also propose  a gradient-based local search (named as Binary Trust Region - BiTR) procedure to solve the binary nonlinear formulation of the assortment optimization problem. There are  major differences between Algorithm \ref{algo:local-search} and the one proposed in \cite{Mai2019Assortment}. First,  our algorithm starts with a greedy heuristic that guarantees an $(1-1/e)$ approximation solution, while there is no performance guarantee for the BiTR. Second, we explore the the structure of the MCP under the GEV family, e.g., the coefficients of the sub-problem's objective function are non-negative and we only care about fixed-size subsets, to build more efficient method to solve the sub-problems.

%%%%%%%%%%%%%%%%%%%%%%%%%%%%%%%
%%%%%%%%%%%%%%%%%%%%%%%%%%%%%%%

\section{Numerical Experiments}
\label{sec:Expr.}
In this section, we provide experimental results to compare our GGX algorithm with existing approaches. We use three datasets from recent literature \citep{Ljubic2018outer, Mai2019Assortment} and generate instances under three popular discrete choice models, i.e., the MNL, MMNL and nested logit models.   
%We prere
%of our binary trust region (BiTR) algorithm on standard data sets from the literature and the comparison between the Gradient-based Local Search approach and some effective approaches mentioned in the previous studies.

\subsection{Experimental Settings}
We will compare our algorithm with the standard greedy heuristic (GH - Step 1 of our GGX algorithm), the multicut  and singlecut outer-approximation algorithms (MOA and OA) \citep{Mai2019Assortment}  and the Branch-and-Cut (BC) \citep{Ljubic2018outer}.
In particular, the for nested logit instances, we only compare the GGX with  GH, OA, MOA approaches, as  BC is not designed to handle such instances.  
Note that for the MNL and MMNL instances, it is possible to formulate the MCP as an MILP and solve by an MILP solver (e.g. IBM's CPLEX). However, as shown in  \cite{Ljubic2018outer} and \cite{Mai2019Assortment}, this approach is  outperformed by the MOA and BC methods. Thus, we do not include the MILP solver in our experiments.

%\textbf{Greedy} We perform only the first phase of the Gradient-based Local Search.
%\textbf{BiTR} We perform phase 1 and 2 of the Gradient-based Local Search method.%
%\textbf{BiTR+LS} We perform the Gradient-based Local Search algorithm.
%\textbf{BC} A branch-and-cut algorithm was proposed by \cite{Ljubic2017} based on a multicut implementation of outer-approximation and sub-modular cuts.
%\textbf{OA} An outer-approximation algorithm was proposed by \cite{Tien2020} based on the cutting plane scheme.
%\textbf{MOA} A multicut outer-approximation was proposed by \cite{Tien2020} that based on the OA algorithm and generate cuts for groups of demand points.
%Our BiTR+LS approach combines two phase such as the Binary Trust Region procedure and the local search procedure (LS). In our implementation, we use the best-improvement schema for the LS procedure, i.e, we consider any move can perform and choose the best one. We denote the set of selected locations as $S$, $S \in M$. This LS procedure includes two operators such as
%\begin{itemize}
%    \item swap(1,1): This operator replaces location $i$, $i \in S$, by location %$j$, $j \notin S$.
%    \item swap(2,2): This operator replaces two locations $i, j$, $(i,j)\in S$, by %two location $k, l$, $(k,l) \notin S$.
%\end{itemize}
%For the computational reason, we perform moves of the swap(1,1) first and just perform moves of the swap(2,2) if there is not any improvement move of swap(1,1).This procedure stops when all possible moves have been tried without improvement. 

We use the following three datasets  as benchmark instances 
and we refer the reader to \cite{Freire2015} for more detailed descriptions. These datasets have been also used in some recent MCP studies \citep{Ljubic2018outer,Mai2019Assortment}.
\begin{itemize}
    \item \textbf{HM14}: The dataset includes 15 problems generated randomly in a plane, with $|I|\in \{50,100,200,400,800\}$ and $m \in \{25,50,100\}$.
    \item \textbf{ORlib}: The dataset includes 11 problems where there are four instances with $(|I|,m)=(50,25)$, four instances with $(|I|,m)=(50,50)$ and three instances with $(|I|,m)=(1000,100)$.
    \item \textbf{PR-NYC} (or \textbf{NYC}): the dataset comes from a large-scale park-and-ride location problem in New York city with $|I| = 82341$ and $m = 59$. As reported in previous studies, these are the largest and most challenging instances.
\end{itemize}
%For the detail of the above data sets, we refer the reader to \cite{Freire2015}.
We employ the same settings of parameters as  in previous studies \citep{Ljubic2017,Tien2020}. The number of facilities that need to be opened $C$ is varied from 2 to 10. The deterministic part of the utility is defined as $v_{ij} = -\beta c_{ij}$ for a location $j\in M$ and $v_{ij^{'}}=-\beta \alpha c_{ij^{'}}$ for each competitor $j^{'}$, where $c_{ij}$ is the distance between zone/client $i\in I$ and location $j\in [m]$, the parameter $\beta$ is the sensitivity of customers about the perceived utilities and $\alpha$ represents the competitiveness of the competitors. These parameters are chosen as $\alpha =\{0.01,0.1,1\}$ and $\beta=\{1,5,10\}$ for datasets HM14 and ORlib, and $\alpha =\{0.5,1,2\}$ and $\beta=\{0.5,1,2\}$ for the NYC dataset. Therefore, for each discrete choice model chosen, each problem above has 81 different instances and the total numbers of instances for HM14, ORlib, NYC are 972, 891, 81, respectively. 

%We should note that there exists some small errors while calculating the objectives because of the large number of digits in the utilities data. Therefore, it is difficult to compare exactly the results from the approaches mentioned above. Thus, we assume that two values are equivalent if the gap between them is less than $0.5\%$. We also set a running time budget of ten minutes for all approaches.  

The experiments are done on a  PC with processor AMD Ryzen 7-3700X CPU @ 3.80 GHz and 16 gigabytes of RAM.  We use MATLAB 2020 to implement and run the algorithms,  and we link to IBM ILOG-CPLEX 12.10 to solve MILPs under default settings. We also take the code used in \cite{Ljubic2017} to generate results for the MNL and MMNL instances with the BC approach.

\subsection{Multinomial Logit - MNL}
We take MNL instances from previous work \citep{Ljubic2018outer,Mai2019Assortment} and report numerical results in Table \ref{table1} below. Each row of the table corresponds to 81 instances and we indicate the largest number of instances solved with the best objective values in bold. We use the same settings as in \cite{Tien2020}. We do not show the CPU times for GH  as it runs very fast. The GH  finishes 26/27 problems  in less than 0.01 seconds and it just needs around 0.5 seconds to finish all  the instances of the largest dataset (i.e. the NYC one). On the other hand, solutions obtained by GH are  relatively good, in the sense that  the percentage gaps between the objective values yielded by those solutions and the best objective values  vary only from 0 to 2.94$\%$. 
%In this experiment, there are only three problems 
%such as "capb", "capc" and "NYC" where BiTR improves the results obtained from Greedy. 
In terms of number instances with the best objective values, GGX performs the best as  it gives the largest number instances with the best objective values in 26/27 problems. Moreover, GGX solves 81/81 instances with the best objective values in 25/27 problems. On the other hand,  GGX only requires short CPU times to finish (the average CPU times are always less than 1.5 seconds except for the NYC instances). Furthermore, when comparing  GGX with the  OA, MOA, and BC approaches, the average CPU times required by GGX are about 78 times lower than  OA, 28 times lower than MOA, and 12 times lower than  BC. For  small instances with $|I|\leq 100$,  BC achieves good performance. It provides the best objective values for all 81 instances of each problem  with the lowest CPU times. However, for larger problem instances ($|I|>100$), BC becomes more expensive, especially for  the three large problems in the ORlib dataset with  $(|I|,|M|)=(800,100)$ and for the NYC instances (the average CPU times are always more than 90 seconds).  The MOA  has the best performance for the NYC dataset, as it only requires 2.32 seconds to give the best objective values for all the 81 instances. In general,  GGX  achieves the best performance for the MNL instances, as compared to the other approaches.

\begin{table}\footnotesize
\centering
\label{table1}
\begin{tabular}{llllllll|llll} 
\hline
\multirow{2}{*}{Instance} & \multicolumn{1}{c}{\multirow{2}{*}{$|I|$ }} & \multicolumn{1}{c}{\multirow{2}{*}{$m$}} & \multicolumn{5}{c|}{\# instances with best objective}    & \multicolumn{4}{c}{Average CPU time (s)}\\ 
\cline{4-12}
   & \multicolumn{1}{c}{} & \multicolumn{1}{c}{}& \multicolumn{1}{l}{GGX} & \multicolumn{1}{l}{GH} & \multicolumn{1}{l}{OA} & \multicolumn{1}{l}{MOA} & \multicolumn{1}{l|}{BC} & \multicolumn{1}{l}{GGX} & \multicolumn{1}{l}{OA} & \multicolumn{1}{l}{MOA} & \multicolumn{1}{l}{BC}  \\ 
\hline
OUR& 50& 25& \textbf{81 }& \textbf{81 } & \textbf{81 } & \textbf{81 }& \textbf{81 } & 0.14 & 19.15    & 0.12 & 0.01 \\ 
\hline
OUR& 50& 50& \textbf{81}& \textbf{81}& 73  & \textbf{81}& \textbf{81}& 0,15 & 109.45   & 0.15 & 0.01 \\ 
\hline
OUR& 50& 100& 79   & 79  & 58  & \textbf{81}& \textbf{81}& 0.24 & 188.91   & 0.32 & 0.05 \\ 
\hline
OUR& 100    & 25& \textbf{81}& 80  & 73  & \textbf{81}& \textbf{81}& 0.14 & 138.30   & 0.18 & 0.01 \\ 
\hline
OUR& 100    & 50& \textbf{81}& 77  & 69  & \textbf{81}& \textbf{81}& 0.15 & 170.80   & 0.28 & 0.03 \\ 
\hline
OUR& 100    & 100& \textbf{81}& 77  & 59  & \textbf{81}& \textbf{81}& 0.26 & 187.43   & 0.60 & 0.13 \\ 
\hline
OUR& 200    & 25& \textbf{81}& \textbf{81}& 72  & \textbf{81}& \textbf{81}& 0.14 & 146.79   & 0.34 & 0.02 \\ 
\hline
OUR& 200    & 50& \textbf{81}& 80  & 64  & \textbf{81}& 80  & 0.16 & 189.46   & 0.80 & 0.06 \\ 
\hline
OUR& 200    & 100& \textbf{81}& 77  & 59  & \textbf{81}& \textbf{81}& 0.33 & 235.27   & 15.99& 0.29 \\ 
\hline
OUR& 400    & 25& \textbf{81}& 73  & 71  & \textbf{81}& 80  & 0.14 & 116.44   & 0,62 & 0.04 \\ 
\hline
OUR& 400    & 50& \textbf{81}& 78  & 60  & 80   & \textbf{81}& 0.18 & 200.32   & 12.13& 0.13 \\ 
\hline
OUR& 400    & 100& \textbf{81}& 73  & 58  & 76   & \textbf{81}& 0.49 & 291.38   & 99.45& 0.65 \\ 
\hline
OUR& 800    & 25& \textbf{81}& 76  & 59  & \textbf{81}& \textbf{81}& 0.14 & 160.71   & 2.27 & 0.11 \\ 
\hline
OUR& 800    & 50& \textbf{81}& 62  & 59  & 64   & 75  & 0.23 & 251.8    & 160.47    & 0.48 \\ 
\hline
OUR& 800    & 100& \textbf{80}& 72  & 55  & 58   & 75  & 0.94 & 363.54   & 234.67    & 14.29\\ 
\hline\hline
cap101 & 50& 25& \textbf{81}   & \textbf{81}  & \textbf{81}  & \textbf{81}   & \textbf{81}  & 0.14 & 0.21& 0.20 & 0.01 \\ 
\hline
cap102 & 50& 25& \textbf{81}   & \textbf{81}  & \textbf{81}  & \textbf{81}   & \textbf{81}  & 0.14 & 0.24& 0.25 & 0.01 \\ 
\hline
cap103 & 50& 25& \textbf{81}   & \textbf{81}  & \textbf{81}  & \textbf{81}   & \textbf{81}  & 0.14 & 0.19& 0.23 & 0.01 \\ 
\hline
cap104 & 50& 25& \textbf{81}& 80  & \textbf{81}& \textbf{81}& \textbf{81}& 0.14 & 0.31& 0.23 & 0.01 \\ 
\hline
cap131 & 50& 50& \textbf{81}& 74  & \textbf{81}& \textbf{81}& \textbf{81}& 0.15 & 0.37& 0.31 & 0.02 \\ 
\hline
cap132 & 50& 50& \textbf{81}   & \textbf{81}  & \textbf{81}  & \textbf{81}   & \textbf{81}  & 0.15 & 0.32& 0.35 & 0.02 \\ 
\hline
cap133 & 50& 50& \textbf{81}& 80  & \textbf{81}& \textbf{81}& \textbf{81}& 0.15 & 0.38& 0.34 & 0.02 \\ 
\hline
cap134 & 50& 50& \textbf{81}   & \textbf{81}  & \textbf{81}  & \textbf{81}   & \textbf{81}  & 0.15 & 0.36& 0.34 & 0.02 \\ 
\hline
capa   & 1000   & 100& \textbf{81}& 58  & \textbf{81}& \textbf{81}& \textbf{81}& 1.31 & 141.10   & 226.68    & 114.74    \\ 
\hline
capb   & 1000   & 100& \textbf{81}& 61  & 68  & 66   & \textbf{81}& 1.23 & 113.27   & 220.45    & 94.72\\ 
\hline
capc   & 1000   & 100& \textbf{81}& 53  & \textbf{81}& \textbf{81}& \textbf{81}& 1.42 & 145.22   & 232.68    & 149.88    \\ 
\hline\hline
NYC& 82341  & 59& \textbf{81}& 72  & 77  & \textbf{81}& 80  & 33.87& 164.01   & 2.32 & 161.71    \\ 
\hline
\hline
Average& \multicolumn{1}{l}{} & \multicolumn{1}{l}{}& 80.89& 75.19    & 71.30    & 78.74& 80.44    & \multicolumn{1}{l}{}    & \multicolumn{1}{l}{}   & \multicolumn{1}{l}{}    & \multicolumn{1}{l}{}    \\
\hline
\end{tabular}
\caption{Numerical results for MNL instances, grouped by the problem name (81 instances per row).}
\end{table}

\subsection{Mixed Logit - MMNL}

In this section, we report numerical results for MMNL instances. 
To generate such instances, we assume that each utility $v_{ij}$, $i\in I, j\in [m]\}$, contains a random error component that follows a normal distribution of zero mean. We also assume that the variance of the random number is proportional to the distance $c_{ij}$. More precisely, each $v_{ij}$ associated with customer zone (or client) $i \in I$ and location $j\in [m]$ is defined  as $v_{ij} = -\theta c_{ij} + c_{ij}\tau_{ij}/3$, where $\tau_{ij}$ is a standard normal random number. We also keep the utilities associated with the competitors deterministic. For each problem, we approximate the objective function by the Monte Carlo method. To do so, we choose a sample size  $N=100$ for the  {HM14} and {ORlib} datasets and $N = 10$ for the NYC one. For the latter, we only choose small $N$ because the NYC problem is already large even under the MNL model.
As mentioned, we consider and solve these MMNL instances as  extended MNL ones, in which the number of customer zones is 5000, 5000 and 823,410 for  instances from the {HM14}, {ORlib} and {NYC} datasets, respectively. We give a CPU time budget  of  600 seconds for  all instances. 
 
 In Table \ref{table2}, for each problem we report the number of instances with the best objective values and the average CPU times over 81 instances  for five approaches, i.e., GGX, GH, OA, MOA, and BC. We indicate in bold the largest numbers of instances solved with the best objective values. The results clearly show that GGX generally outperforms other approaches. More precisely, GGX  manages to return best objective values for all instances considered (i.e. 2187/2187 instances). Moreover,  GH also performs  very well, in the sense that the percentage gaps between the objective values given by  GH and the best objective values only vary from 0\% to 2.92$\%$.  
 %The number instances with the best objective of Greedy, BiTR, and BC is quite similar when there are 11/27, 11/27, and 12/27 problems respectively, where the number instances with the best objective is maximum (i.e, 81/81 instances with the best objective on each problem) and the average number instances with the best objective are 75.81, 75.85, and 77.56 respectively. 
 In terms of CPU time, GH  is still the fastest approach when it just requires less than 2 seconds to solve every instance except the NYC ones, which take only about 6.5 seconds in average. The GGX approach, even though being slower than the GH, but is still much faster than the others.  
 We also observe that the OA, MOA, and BC approaches need much more time to solve MMNL instances, as compared to solving the MNL instances. The average CPU times required by these three approaches are more than 250 seconds. 
 In overall, GGX dominates GH  in terms of returned objective value, and outperforms
 OA, MOA and BC  in terms of both returned objective value and CPU time.

\begin{table}\footnotesize
\centering
\begin{tabular}{llllllll|lllll} 
\hline
\multirow{2}{*}{Problem} & \multirow{2}{*}{$|I|$ } & \multirow{2}{*}{$m$ } & \multicolumn{5}{c|}{\# instances with best objective}& \multicolumn{5}{c}{Average CPU time (s)}\\ 
\cline{4-13}
 & & & GGX   & GH & OA    & MOA   & BC & GGX    & GH   & OA& MOA    & BC \\ 
\hline
OUR   & 50   & 25   & \textbf{81}& \textbf{81}  & 66    & 58    & \textbf{81}  & 0.24   & 0.01 & 236.36 & 189.23 & 1.66    \\ 
\hline
OUR   & 50   & 50   & \textbf{81}& \textbf{81}  & 46    & 54    & \textbf{81}  & 0.83   & 0.02 & 312.73 & 232.16 & 6.45    \\ 
\hline
OUR   & 50   & 100  & \textbf{81}& \textbf{81}  & 41    & 29    & \textbf{81}  & 3.49   & 0.05 & 390.60 & 425.57 & 40.78   \\ 
\hline
OUR   & 100  & 25   & \textbf{81}& 75 & 53    & 55    & \textbf{81}  & 0.29   & 0.01 & 279.68 & 212.33 & 11.37   \\ 
\hline
OUR   & 100  & 50   & \textbf{81}& \textbf{81}  & 47    & 31    & \textbf{81}  & 1.11   & 0.03 & 356.28 & 430.85 & 36.57   \\ 
\hline
OUR   & 100  & 100  & \textbf{81}& 75 & 41    & 22    & 78 & 4.76   & 0.10 & 528.64 & 519.79 & 133.67  \\ 
\hline
OUR   & 200  & 25   & \textbf{81}& \textbf{81}  & 54    & 33    & \textbf{81}  & 0.35   & 0.02 & 317.28 & 375.46 & 47.32   \\ 
\hline
OUR   & 200  & 50   & \textbf{81}& 77 & 39    & 27    & 78 & 1.51   & 0.05 & 402.01 & 464.52 & 126.53  \\ 
\hline
OUR   & 200  & 100  & \textbf{81}& 77 & 34    & 28    & 77 & 6.75   & 0.23 & 456.24 & 328.89 & 174.75  \\ 
\hline
OUR   & 400  & 25   & \textbf{81}& 77 & 52    & 32    & 77 & 0.44   & 0.03 & 336.85 & 404.77 & 138.9   \\ 
\hline
OUR   & 400  & 50   & \textbf{81}& 74 & 40    & 25    & 74 & 2.03   & 0.13 & 472.38 & 515.90 & 237.10  \\ 
\hline
OUR   & 400  & 100  & \textbf{81}& 63 & 29    & 15    & 63 & 9.34   & 0.67 & 569.4  & 570.10 & 325.16  \\ 
\hline
OUR   & 800  & 25   & \textbf{81}& 66 & 52    & 29    & 68 & 0.60   & 0.07 & 345.27 & 396.63 & 229.83  \\ 
\hline
OUR   & 800  & 50   & \textbf{81}& 73 & 39    & 27    & 74 & 3.06   & 0.34 & 406.18 & 467.89 & 341.22  \\ 
\hline
OUR   & 800  & 100  & \textbf{81}& 72 & 32    & 18    & 72 & 14.29  & 1.25 & 562.63 & 565.37 & 421.74  \\ 
\hline\hline
cap101& 50   & 25   & \textbf{81}& 70 & 57    & 54    & 77 & 0.24   & 0.01 & 388.27 & 485.02 & 276.30  \\ 
\hline
cap102& 50   & 25   & \textbf{81}& 72 & 57    & 54    & 80 & 0.24   & 0.01 & 388.19 & 497.24 & 268.27  \\ 
\hline
cap103& 50   & 25   & \textbf{81}& 66 & 50    & 51    & 78 & 0.24   & 0.01 & 365.41 & 457.82 & 218.40  \\ 
\hline
cap104& 50   & 25   & \textbf{81}& 67 & 56    & 51    & 74 & 0.25   & 0.01 & 369.85 & 484.69 & 268.10  \\ 
\hline
cap131& 50   & 50   & \textbf{81}& \textbf{81}  & 33    & 38    & \textbf{81}  & 0.84   & 0.02 & 459.78 & 619.15 & 425.71  \\ 
\hline
cap132& 50   & 50   & \textbf{81}& 80 & 33    & 38    & 80 & 0.84   & 0.02 & 452.13 & 620.70 & 422.28  \\ 
\hline
cap133& 50   & 50   & \textbf{81}& \textbf{81}  & 33    & 35    & \textbf{81}  & 0.83   & 0.02 & 459.68 & 643.83 & 435.96  \\ 
\hline
cap134& 50   & 50   & \textbf{81}& \textbf{81}  & 33    & 40    & \textbf{81}  & 0.84   & 0.02 & 458.49 & 652.69 & 429.13  \\ 
\hline
capa  & 1000 & 100  & \textbf{81}& \textbf{81}  & 14    & 14    & \textbf{81}  & 16.50  & 1.52 & 677.73 & 677.73 & 600.00  \\ 
\hline
capb  & 1000 & 100  & \textbf{81}& \textbf{81}  & 44    & 14    & \textbf{81}  & 16.46  & 1.52 & 600.91 & 691.05 & 600.00  \\ 
\hline
capc  & 1000 & 100  & \textbf{81}& \textbf{81}  & 38    & 13    & \textbf{81}  & 16.47  & 1.52 & 600.91 & 637.31 & 600.00  \\ 
\hline\hline
NYC   & 82341& 59   & \textbf{81}& 72 & 75    & 68    & 72 & 114.42 & 6.50 & 178.18 & 135.24 & 381.48  \\ 
\hline
\hline
Average   & & & \textbf{81}  & 75.81   & 44.00 & 35.30 & 77.56   &   & &   &   &    \\
\hline
\end{tabular}
\caption{Numerical results for MMNL instances, grouped by the problem name (81 instances per row).}
\label{table2}
\end{table}

\subsection{Nested Logit Model}
This section reports numerical results for nested logit instances. We perform a comparison between 4 approaches, namely, the GGX, GH, OA and MOA algorithms.
We do not include the BC approach in this experiment, as it is not designed to handle nested logit instances. For the OA and MOA approaches, since it is quite straightforward to generate outer-approximation cuts using gradient information, we apply these algorithms to solve the nested logit instances to see how they perform. Note that, in the context, the objective function is no-longer concave, thus OA and MOA become heuristic with no performance guarantee, to the best of our knowledge.     
To generate nested logit instances, we build a customer nested logit model by 
partitioning the set of locations into $L=5$ different and disjoint groups of equal size. 
%We denote $g_{n}$ as the $n$-th group with $n=1,2,...,5$. For example, we assume that the number of locations $|M|$ is 10. Then, we have 5 groups as $g_1=\{1,2\}$; $g_2=\{3,4\}$; $g_3=\{5,6\}$; $g_4=\{7,8\}$; $g_5=\{9,10\}$. 
In particular, the NYC dataset has 59 locations ($m=59$), so for this problem we partition the locations into four groups with 10 locations and one with 9  locations. We also choose the nested logit parameters as $\mu = (1.1,1.2,1.3,1.4,1.5)$, noting that more nests and/or other nested logit parameters can be chosen. Our selections here are just to illustrate the performance of different algorithms in handling  GEV instances.  We also give a time budget of 600 seconds for all the algorithms. 

Table \ref{table3} reports comparison results of the four approaches. Each row of the table corresponds to 81 solved instances and we also indicate the largest numbers of instances solved with best objective values in bold. The results clearly show that GGX  outperforms the other  approaches in terms of the number of instances solved with the best objective values. More precisely, GGX gives the best objective values for all problem instances while GH only performs the best for 9/27 problems. In terms of CPU time,  GGX  is not very fast. In particular, for the NYC instances, the average CPU time is about 355.36 seconds and is much larger than the average CPU times required by GH, OA and MOA. The reason is that the objective function in this context is quite expensive to evaluate, as compared to
the cases of the  MMNL and MNL models, and  the \textit{exchanging} procedure of the GGX (Phase 3) requires calculating the objective function several times to find a pair of locations to swap. The GH is still very fast and the returned objective values are pretty close to the best values given by GGX. The percentage gaps between the objective values obtained from GH and the best objective values only vary from 0 to 3.32$\%$. 
%Therefore, BiTR improves the results obtained from Greedy in only one problems, the "cap133" problem.
The OA and MOA approaches, even though run very fast, but give bad solutions. This can be explained by the fact that the objective function under a nested logit model is highly non-concave, thus a subgradient cut (or an outer-approximation cut) could potentially remove good solutions during the cutting-plane procedure.

\begin{table}\footnotesize
\centering
\begin{tabular}{lllllll|llll} 
\hline
\multicolumn{1}{c}{\multirow{2}{*}{Problem}} & \multicolumn{1}{c}{\multirow{2}{*}{$|I|$ }} & \multicolumn{1}{c}{\multirow{2}{*}{$m$ }} & \multicolumn{4}{c}{\# instances with best objective}& \multicolumn{4}{c}{Average CPU time (s)}   \\ 
\cline{4-11}
\multicolumn{1}{c}{}      & \multicolumn{1}{c}{}     & \multicolumn{1}{c}{}     & \multicolumn{1}{l}{GGX}             & \multicolumn{1}{l}{GH} & \multicolumn{1}{l}{OA} & \multicolumn{1}{l|}{MOA} & \multicolumn{1}{l}{GGX} & \multicolumn{1}{l}{GH} & \multicolumn{1}{l}{OA} & \multicolumn{1}{l}{MOA}  \\ 
\hline
OUR& 50& 25& \textbf{81}      & \textbf{81}            & 24  & 20    & 0.35 & 0.01& 0.15& 0.09  \\ 
\hline
OUR& 50& 50& \textbf{81}      & \textbf{81}            & 2   & 0     & 1.94 & 0.03& 2.83& 0.12  \\ 
\hline
OUR& 50& 100   & \textbf{81}      & \textbf{81}            & 27  & 27    & 13.08& 0.10& 104.39                 & 0.19  \\ 
\hline
OUR& 100   & 25& \textbf{81}      & 80  & 18  & 15    & 0.43 & 0.01& 0.66& 0.16  \\ 
\hline
OUR& 100   & 50& \textbf{81}      & 75  & 3   & 0     & 2.44 & 0.04& 0.39& 0.16  \\ 
\hline
OUR& 100   & 100   & \textbf{81}      & 76  & 0   & 1     & 16.04& 0.11& 83.31                  & 0.27  \\ 
\hline
OUR& 200   & 25& \textbf{81}      & \textbf{81}            & 2   & 0     & 0.57 & 0.02& 0.80& 0.18  \\ 
\hline
OUR& 200   & 50& \textbf{81}      & 80  & 1   & 0     & 3.41 & 0.19& 1.85& 0.23  \\ 
\hline
OUR& 200   & 100   & \textbf{81}      & \textbf{81}            & 0   & 0     & 20.76& 0.15& 111.70                 & 0.46  \\ 
\hline
OUR& 400   & 25& \textbf{81}      & 72  & 24  & 20    & 0.85 & 0.16& 0.75& 0.18  \\ 
\hline
OUR& 400   & 50& \textbf{81}      & 77  & 5   & 0     & 5.37 & 0.21& 1.83& 0.35  \\ 
\hline
OUR& 400   & 100   & \textbf{81}      & 72  & 0   & 0     & 32.02& 0.08& 167.55                 & 1.15  \\ 
\hline
OUR& 800   & 25& \textbf{81}      & 77  & 0   & 0     & 0.93 & 0.04& 0.89& 0.31  \\ 
\hline
OUR& 800   & 50& \textbf{81}      & 66  & 0   & 0     & 6.38 & 0.10& 2.14& 0.61  \\ 
\hline
OUR& 800   & 100   & \textbf{81}      & 69  & 22  & 18    & 41.99& 0.30& 299.35                 & 5.64  \\ 
\hline
cap101 & 50& 25& \textbf{81}      & \textbf{81}            & 1   & 0     & 0.35 & 0.01& 0.09& 0.10  \\ 
\hline
cap102 & 50& 25& \textbf{81}      & \textbf{81}            & 1   & 0     & 0.35 & 0.01& 0.10& 0.10  \\ 
\hline
cap103 & 50& 25& \textbf{81}      & \textbf{81}            & 0   & 0     & 0.35 & 0.01& 0.09& 0.10  \\ 
\hline
cap104 & 50& 25& \textbf{81}      & 78  & 0   & 0     & 0.35 & 0.01& 0.09& 0.09  \\ 
\hline
cap131 & 50& 50& \textbf{81}      & 74  & 0   & 0     & 1.96 & 0.03& 0.10& 0.11  \\ 
\hline
cap132 & 50& 50& \textbf{81}      & \textbf{81}            & 0   & 0     & 1.96 & 0.03& 0.11& 0.11  \\ 
\hline
cap133 & 50& 50& \textbf{81}      & 78  & 0   & 0     & 1.97 & 0.03& 0.11& 0.11  \\ 
\hline
cap134 & 50& 50& \textbf{81}      & 80  & 0   & 0     & 1.96 & 0.03& 0.10& 0.11  \\ 
\hline
capa   & 1000  & 100   & \textbf{81}      & 60  & 0   & 0     & 46.17& 0.37& 1.02& 66.50 \\ 
\hline
capb   & 1000  & 100   & \textbf{81}      & 54  & 0   & 0     & 46.23& 0.37& 10.91                  & 66.24 \\ 
\hline
capc   & 1000  & 100   & \textbf{81}      & 50  & 0   & 0     & 46.18& 0.37& 8.38& 56.44 \\ 
\hline
NYC& 82341 & 59& \textbf{81}      & 72  & 5   & 4     & 355.36                  & 6.89& 0.85& 0.65  \\ 
\hline
\multicolumn{1}{c}{Average}                  & \multicolumn{1}{c}{}     & \multicolumn{1}{c}{}     & \multicolumn{1}{l}{\textbf{81} } & 74.78                  & 5.00& \multicolumn{1}{l}{3.89} &      &     &     &       \\
\hline
\end{tabular}
\caption{Numerical results for nested-logit instances, grouped by the problem name (81 instances per row).}
\label{table3}
\end{table}

We look more closely to the NYC problem (the largest problem) to see how the algorithm works. In table \ref{table4}, we report comparison results for the NYC instances in detail. Each row of the table corresponds to 9 instances with a value of $C$, varying from 2 to 10.  GGX performs the best in terms of objective value, as it gives best objective values for all the instances while GH only gives 6/9 best objective values for $C \in \{2,3,4\}$. The numbers of instances with best objective values given by OA and MOA  are  very low. They both have 4 instances with the best objective values when $C=2$ and the OA has one more instance with the best objective value when $C=5$. This clearly shows that OA and MOA are outperformed by GGX and GH. On the other hand, in terms of CPU time, GGX is much more expensive than the other approaches. The average CPU times required by  GGX  is about 52 times, 418 times, and 547 times higher than those required by  GH, OA, and MOA approaches, respectively.
In summary, for these large instances,  GH performs much better as compared to  OA and MOA, and GGX  manages to significantly improve the objective values returned by  GH.  

\begin{table}
\centering
\begin{tabular}{c|cccc|cccc} 
\hline
\multicolumn{1}{c|}{\multirow{2}{*}{C}} & \multicolumn{4}{c|}{\begin{tabular}[c]{@{}c@{}}\# instances \\with best objective values\end{tabular}}   & \multicolumn{4}{c}{\begin{tabular}[c]{@{}c@{}}Average\\~CPU time (s) \end{tabular}}\\ 
\cline{2-9}
\multicolumn{1}{c|}{} & \multicolumn{1}{c}{GGX} & \multicolumn{1}{c}{GH} & \multicolumn{1}{c}{OA} & \multicolumn{1}{c|}{MOA} & \multicolumn{1}{c}{GGX} & \multicolumn{1}{c}{GH} & \multicolumn{1}{c}{OA} & \multicolumn{1}{c}{MOA}  \\ 
\hline
2& \textbf{9}& 6 & 4 & 4   & 81.90 & 2.28 & 1.03 & 0.74\\ 
\hline
3& \textbf{9}& 6 & 0 & 0   & 234.77& 3.53 & 1.24 & 0.83\\ 
\hline
4& \textbf{9}& 6 & 0 & 0   & 454.85& 4.75 & 1.26 & 0.94\\ 
\hline
5& \textbf{9}& \textbf{9}& 1 & 0   & 551.34& 5.74 & 1.13 & 0.67\\ 
\hline
6& \textbf{9}& \textbf{9}& 0 & 0   & 556.38& 6.94 & 1.14 & 0.60\\ 
\hline
7& \textbf{9}& \textbf{9}& 0 & 0   & 541.49& 8.04 & 0.52 & 0.54\\ 
\hline
8& \textbf{9}& \textbf{9}& 0 & 0   & 538.21& 9.19 & 0.49 & 0.52\\ 
\hline
9& \textbf{9}& \textbf{9}& 0 & 0   & 537.17& 10.27& 0.48 & 0.51\\ 
\hline
10   & \textbf{9}& \textbf{9}& 0 & 0   & 541.68& 11.25& 0.41 & 0.52\\ 
\hline
\multicolumn{1}{l|}{Average}& \textbf{9}& 8 & 0.56 & 0.44& \multicolumn{1}{l}{}    & \multicolumn{1}{l}{}   & \multicolumn{1}{l}{}   & \multicolumn{1}{l}{}\\
\hline
\end{tabular}
\caption{Comparison results for NYC instances, grouped by $C$, 9 instances per row.}
\label{table4}
\end{table}

\section{Conclusion}
\label{sec:conc}

In this paper we have studied the maximum capture problem in facility location where customer behavior is captured by any GEV model. By leveraging the properties of the GEV generating function, we have showed that the objective function is monotonic   and submodular, implying that a simple greedy heuristic can always give a solution whose value is at least $(1-1/e)$ times the optimal values. We have further developed an algorithm based on a greedy heuristic, a gradient-based local search and an exchanging procedure to solve the problem under any GEV model and the MMNL model. We have tested and compared our algorithm with some state-of-the-art algorithms using  MNL, MMNL and nested logit instances and our numerical experiments clearly demonstrate the advantages of our approach, in terms of both returned objective value and CPU time.  

Our theoretical findings and algorithm can be applied to the maximum capture problem under any GEV model, including the popular MNL model and other complex GEV models in the literature. Future directions would be to formulate and solve a maximum capture  problem in the situation that the  choice parameters are not known with certainty, or to consider a combination of facility location and security planning under the MNL/MMNL or any GEV models.  

% Acknowledgments here
%\ACKNOWLEDGMENT{...}

\bibliographystyle{plainnat}
\bibliography{reference,refs}

\end{document}